\documentclass[sn-mathphys]{sn-jnl}
\jyear{2021}%
\def\pput(#1,#2)#3{\noindent\smash{\raise#2pt\hbox to 0pt
   {\kern #1pt #3\hss}}\ignorespaces}
\usepackage{graphicx,upquote,mleftright}
\def\Im{\hbox{\rm Im\kern .7pt}}
\def\Re{\hbox{\rm Re\kern .7pt}}
\def\dO{\partial\kern .3pt\Omega}
\def\dG{\partial\kern .3pt G}\def\dQ{\partial\kern .3pt Q}
\def\zb{{\bf 0}}
\def\Qm{Q_-}
\def\G{{\bf G}}\def\Q{{\bf Q}}
\def\npoles{n_{\hbox{\tiny poles}}}
\def\npoly{n_{\hbox{\tiny poly}}}
\def\matlab{\hbox{MATLAB}}

\title[Rectangular eigenproblems]{Rectangular eigenvalue problems}

\author[1]{\fnm{Behnam} \sur{Hashemi}}\email{bhashemi@qu.edu.qa}
\author[2]{\fnm{Yuji} \sur{Nakatsukasa}}\email{nakatsukasa@maths.ox.ac.uk}
\author[2]{\fnm{Lloyd N.} \sur{Trefethen}}\email{trefethen@maths.ox.ac.uk}
\affil[1]{\orgdiv{Mathematics Program, Dept.\ of Mathematics, Statistics and Physics, College
of Arts and Sciences},
\orgname{Qatar University}, \orgaddress{2713}, \city{Doha}, \country{Qatar}}
\affil[2]{\orgdiv{Mathematical Institute}, \orgname{University of Oxford},
\orgaddress{\street{Woodstock Rd.}, \city{Oxford}, 
\postcode{OX2 6GG}, \country{UK}}}

\usepackage{graphicx,upquote}

\begin{document}

\abstract{Often the easiest way to discretize an ordinary or partial
differential equation is by a {\em rectangular numerical method,}
in which $n$ basis functions are sampled at $m \gg n$ collocation points.
We show how eigenvalue problems can be solved
in this setting by QR reduction to square
matrix generalized eigenvalue problems.  The method applies equally in the
limit ``\kern .7pt $m=\infty$'' of eigenvalue problems for quasimatrices.
Numerical examples are presented as
well as pointers to some related literature.}



\maketitle 


\section{\label{introd}Introduction}
Problems involving ordinary and partial differential equations
(ODE\kern .5pt s and PDE\kern .5pt s) are traditionally discretized by
square matrices.  Such methods are effective when a well-conditioned
basis is available in which to expand the numerical solution and
good quadrature or collocation points are known at which to enforce
the equations.  Sometimes, however, these conditions do not hold,
and it becomes advantageous to sample the equation at more data
points than there are basis functions and to solve the problem in a
least-squares formulation.  We call these {\em rectangular numerical
methods.} The aim of this paper is to propose rectangular numerical
methods for ODE and PDE eigenvalue problems.

Rectangular numerical methods have appeared in many areas, though
they have rarely taken center stage.  When Fourier, Chebyshev, or related
expansions are involved, one can speak of rectangular spectral
methods~\cite{drischale}, though Boyd observes that such methods
are ``relatively uncommon''~\cite[sec.~3.1]{boyd}.  In the finite
elements literature there are Least-Squares Finite
Element Methods~\cite{lsfem,jiang,monk}.  With expansion functions that satisfy the
differential equation but not the boundary conditions, one gets
series methods~\cite{series} or the Method of Fundamental Solutions
(MFS)~\cite{bb,fk} or lightning or log-lightning methods for PDE
problems with corner singularities~\cite{lightning,loglightning}.
Related expansions that do not satisfy the differential equation
and hence
need fitting in the interior of a domain, not just on the boundary,
lead to least-squares methods for radial basis functions
(RBFs) or other kernels~\cite{buhmann,fasshauer,km,pd06}.  RBF methods are an
example of the broad category of {\em meshfree} methods.

Our plan is to set forth some of the simplest methods for
solving rectangular eigenproblems and illustrate them with a sequence of
examples.  The closest previous contributions we know of on this
topic are by Manzhos and coauthors, who have developed what they call ``rectangular
collocation'' methods for eigenvalue problems in quantum chemistry~\cite{km,myc}, and
by the first two authors~\cite{hn}.  The emphasis in~\cite{hn}
is on spectral methods for ODE\kern .5pt s, and the linear algebra is
carried out by the method of Ito and Murota~\cite{im}, involving the
singular value decomposition (SVD) of a matrix with twice as many
columns as there are basis functions.  (Important earlier related
papers are~\cite{boutry} and~\cite{wright}.)  Here we look at a wider
range of problems and propose simpler methods of linear algebra based
on the QR decomposition of a matrix without the doubled dimension.

We will mainly deal with fully discrete $m\times n$ rectangular
matrices, always with $m>n$.  As pointed out in~\cite{hn}, however,
it makes good sense conceptually to consider the limit in which
the columns are functions of one or more continuous variables, so
that instead of matrices, we have quasimatrices; see~\cite{bt,house}
and~\cite[chap.~6]{chebfun}.  Nothing essential changes here, so we
shall simply include quasimatrices in the discussion as the case
``$m=\infty$''.  For spectral ODE problems, the quasimatrices can
be realized numerically in Chebfun~\cite{chebfun}, and the first two
of the examples of Section 3 follow this path.  After that, our computed
examples are fully discrete, though the mathematical derivations
apply equally to $m<\infty$ or $m=\infty$.

Rectangular numerical methods for eigenvalue problems are related
to ideas going back a century, first associated with Rayleigh, Ritz, and 
Galerkin, in
which square matrix approximations are obtained by quadrature and
projection~\cite{saad}.  (A fasinating historical discussion is given
in~\cite{ganderwanner}.)
In the finite elements literature, Galerkin and
Petrov--Galerkin methods can often be interpreted this way.
Arnoldi and Jacobi-Davidson iterative
methods for computing eigenvalues of large matrices are also of this
nature.  What is different in the present paper is that no explicit quadrature or
projection ideas are employed, just numerical algorithms applied
to eigenvalue problems configured rectangularly.  This diminishes the
need for case-by-case analysis and permits
great flexibility in the choice of basis functions and sample points.

\section{\label{method}The numerical method, three variants}
Let $L$ be a linear operator acting on functions 
in a univariate or multivariate domain $\Omega$,
and suppose we seek eigenvalues $\lambda$ and nonzero eigenfunctions
$u$ such that
\begin{equation}
Lu = \lambda u. 
\label{diffeq}
\end{equation}
We shall consider three variants of this problem, in which (\ref{diffeq})
is coupled with no boundary conditions, a finite number of boundary conditions,
or boundary conditions applied on a continuum.
In all three cases we suppose that for
some $n\ge 1$, we have a set of functions $g_1, \dots, g_n$ defined in $\Omega$
whose span contains good approximations to the eigenfunctions of
interest, and we let
$G$ be the $m\times n$ matrix whose columns are these functions.
(If $m<\infty$, the columns consist of samples of the functions at $m$ points.)
Setting 
\begin{equation}
u = Gx,
\label{ansatz}
\end{equation}
we seek a coefficient vector $x\in \mathbb{C}^n$ such that
\begin{equation}
LGx = \lambda \kern .5pt Gx. 
\label{Geq}
\end{equation}
This is an $m\times n$ generalized eigenvalue problem, which can
also be described as the eigenvalue problem for the $m\times n$
matrix or quasimatrix pencil $LG-\lambda \kern .5pt G$.  Like most
rectangular eigenproblems, it will not have exact solutions in
general~\cite{boutry,im,wright}, but if $G$ is well chosen, we expect it
to have solutions accurate to many digits of accuracy.
In such contexts $G$ will often be highly ill-conditioned.

Note that if $m<\infty$, so that the columns of $G$ are discrete
samples of the functions $g_k$ and $u$ is also a discrete vector,
then the expression $LG$ in (\ref{Geq}) cannot be exactly the
product of $L$ and $G$.  Instead it should be interpreted as the
discrete matrix whose columns come from sampling the functions
$Lg_k$ at $m$ points in $\Omega$.  The same qualification applies to the product $BG$
introduced in (\ref{bc2}) below.

\smallskip
{\em Variant 1.~No boundary conditions.}
Suppose first that (\ref{diffeq}) is the whole problem: 
as in~\cite{boutry} and~\cite{im}, there are no explicit boundary conditions.
An example would be the\label{p3}
harmonic oscillator $-u'' + x^2 u = \lambda u$ defined on the
real axis, with eigenvalues $1,3,5,\dots.$  
In this case our proposed numerical method begins by
computing the QR factorization of $G$,
\begin{equation}
G = QR,
\label{qr}
\end{equation}
where $Q$ is $m\times n$ and $R$ is $n\times n$ and upper-triangular.
(For details of QR factorization in
the quasimatrix case $m=\infty$, see~\cite{house}.)
Premultiplying (\ref{Geq}) by $Q^*$ then gives
\begin{equation}
Q^*(LG\kern .7pt)x = \lambda \kern .5pt Q^*\kern -.7pt Gx = \lambda Rx.
\label{new}
\end{equation}
This equation enforces the condition that the
residual $(LG\kern .7pt)x - \lambda \kern .5pt Gx$ is orthogonal to the
range of $G$.
(Alternatively one could require
\begin{equation}
G^*(LG\kern .7pt)x = \lambda \kern .5pt G^*\kern -.7pt Gx,
\label{unstable}
\end{equation}
though with less numerical stability, in analogy to the normal
equations for least-squares fitting problems; compare~\cite[eq.~(6)]{km} and the
earlier~\cite[eq.~(5)]{nakatsuji}.)
Equation~(\ref{new}) is a square matrix generalized eigenvalue problem with
dimensions $n\times n$, which we solve by the standard QZ algorithm.
Note that (\ref{Geq}) implies (\ref{new}).  Conversely,
(\ref{new}) implies (\ref{Geq}) if
the columns of $LG$ lie in the column space of $G$.  This may
or may not hold exactly, but in many applications it will hold to
high accuracy, making (\ref{Geq}) and (\ref{new}) effectively equivalent.

In the computation above, as in Variants~2 and~3 below, an alternative
(mathematically equivalent) possibility
is to use the SVD instead of the QR decomposition to construct
an orthonormal basis of the columns of $G$.  In our experience this
may improve the accuracy slightly, typically by less than one digit, at the cost
of a slight increase in computing time.  
We have not investigated the matter carefully.

\smallskip
{\em Variant 2.~Finite set of boundary conditions.}
Suppose next that (\ref{diffeq}) is coupled with a finite set\label{p3b}
of $\mu > 0$ homogeneous linear boundary conditions, as is
considered (along with other possibilities) in~\cite{hn}.
An example would be
$-(4/\pi^2)u'' = \lambda u$ on $[-1,1$] with boundary conditions $u(\pm 1) = 0$, with
eigenvalues $1,4,9,\dots.$
We can write the boundary conditions in the form
\begin{equation}
Bu = \zb,
\label{bc1}
\end{equation}
where $B$ (``boundary'') is a $\mu\times m$ matrix or row quasimatrix and
$\zb$ is the $\mu\times 1$ zero vector.
(In the quasimatrix
case each row of $B$ is a linear functional, which might, for example, evaluate
$u$ or $u'$ at a boundary point.)
Applying (\ref{ansatz}), this becomes the  
$\mu\times n$ set of equations
\begin{equation}
(BG\kern .7pt)x = \zb.
\label{bc2}
\end{equation}
We can now combine (\ref{bc2}) with (\ref{Geq})
to get the $(m+\mu)\times n$ rectangular generalized eigenvalue problem
\begin{equation}
\mleft[
\begin{array}{c}
 \\ LG \\ \\ \hline \\[-10pt] BG
\end{array}
\mright]
x
= \lambda
\mleft[
\begin{array}{c}
 \\ G \\ \hphantom{LG}\\ \hline \\[-10pt] 0
\end{array}
\mright]
x,
\label{var2}
\end{equation}
where $0$ denotes the zero matrix of dimensions $\mu\times n$.
Various methods can be employed to make this equation square, as
discussed in~\cite{hn} in the context of the Ito-Murota formulation.  The simplest,
analogous to what is called the ``tau method'' of imposing boundary
conditions in spectral methods~\cite{boyd}, is to let $\Qm$ denote the $m\times (n-\mu)$
matrix or quasimatrix consisting of $Q$ with its final $\mu$ columns removed and then
consider 
\begin{equation}
\mleft[
\begin{array}{c}
 \\ \Qm^*LG \\ \\ \hline \\[-10pt] BG
\end{array}
\mright]
x
= \lambda
\mleft[
\begin{array}{c}
 \\ \Qm^*G \\ \hphantom{\Qm^*LG}\\ \hline \\[-10pt] 0
\end{array}
\mright]
x.
\label{var2again}
\end{equation}
This equation enforces the boundary conditions exactly while
requiring the residual 
$(LG\kern .7pt)x - \lambda \kern .5pt Gx$ to be orthogonal
to the range of the first $n-\mu$ columns of $G$.  It is a
square matrix generalized eigenvalue problem of
dimensions $n\times n$, which again we solve by standard methods.
For essentially the same structure but not based on a QR
factorization, see~\cite[sec.~5]{drischale} and~\cite[sec.~5]{aurentz}.

\smallskip
{\em Variant 3.~Continuum of boundary conditions.}
Finally, suppose (\ref{diffeq}) is coupled with a continuum of
homogeneous linear boundary conditions.  Specifically, suppose  we
have a PDE in a domain $\Omega$ of dimension $d\ge 2$ and a boundary
condition applied on the boundary $\dO$ of dimension $d-1$.  An example would be $-\Delta
u = \lambda u$ on the unit disk with boundary condition $u=0$ on
the unit circle, whose first eigenvalue is $5.7831859629\dots$,
the square of the smallest root of the Bessel function $J_0(x)$.
\label{bessel}

In this continuous case equations (\ref{bc1})--(\ref{var2})
continue to apply, but the meaning of the $\mu\times m$
boundary condition matrix $B$ is adjusted.  Now $\mu$ will be
either $\infty$, if we formulate the boundary conditions in a
continuous manner, or a large finite number, if we discretize.
Upon multiplying by $G$ we get an object $BG$ of dimensions
$\mu\times n$.  Now that $\mu$ is large or infinite, it
is no longer appropriate to attempt to enforce the boundary
conditions exactly.  Instead, the natural thing to do is to treat
all of (\ref{var2}), both the upper and lower parts, in a
least-squares fashion.  

We do this as follows.  Let $\G$ denote the $(m+\mu)\times n$ structure
\begin{equation}
\G = 
\mleft[
\begin{array}{c}
 \\ G \\ \hphantom{LG}\\ \hline \\[-10pt] \dG
\end{array}
\mright],
\label{G}
\end{equation}
whose columns below the line correspond to $G$ restricted
to the boundary---the boundary traces of the columns of $G$.  Thus each
column of $\G$ is an object whose upper part is a function of $d$ dimensions (or
its $m$-point discretization if $m<\infty$) and whose lower part is a function of $d-1$ dimensions
(or its $\mu$-point discretization if $\mu<\infty$).  We now compute a QR factorization
of $\G$,
\begin{equation}
\G = \Q R,
\label{qrG}
\end{equation}
in which $\Q$ has the same structure as $\G$,
\begin{equation}
\Q = 
\mleft[
\begin{array}{c}
 \\ Q \\ \hphantom{LG}\\ \hline \\[-10pt] \dQ
\end{array}
\mright],
\label{Q}
\end{equation}
and $R$ is an $n\times n$ upper-triangular matrix.
The columns of $\Q$ are orthonormal with respect to an inner 
product that combines integrals (or sums, when $m$ or $\mu$ is finite) associated with
both $\Omega$ and $\dO$.  For discussions of QR factorization and singular
value decomposition of such
mixed objects, see~\cite{hn}.   One could analyze what relative
weighting is most appropriate in balancing the two halves,
but the expectation is that in applications
it will not make much difference.  In our fully discrete computed
examples we give equal weights to all sample points, both the
$m$ points in the interior and the $\mu$ points on the boundary.

To square up the eigenvalue problem, we left-multiply (\ref{var2})
by $\Q^*$ to obtain
\begin{equation}
\left[ Q^*\kern -1pt  (LG\kern .7pt) + (\dQ)^*\kern -1pt  (BG\kern .7pt)\kern .7pt \right] x = 
\lambda \kern .5pt Q^*\kern -1pt  Gx .
\label{var3}
\end{equation}
This equation enforces the condition that a combined boundary-interior residual
is orthogonal in the mixed inner product to the basis vectors (columns of $G$\kern .5pt) and their
boundary traces (columns of $\dG$\kern .5pt).  Like (\ref{new}) and (\ref{var2again}),
(\ref{var3}) is an $n\times n$ generalized eigenvalue problem, and again we solve it
by standard numerical methods.

We now turn to computed examples.  The ODE problems of section~\ref{proto} 
illustrate variants 1 and 2, and the PDE problems of section~\ref{2D} illustrate
variant~3.

\section{\label{proto}One-dimensional examples (ODE\kern .5pt s)}

~\par \vskip -14pt

\indent{\em Example 1.  Harmonic oscillator with no boundary conditions.}
We begin with the harmonic oscillator mentioned on p.~\pageref{p3},
$-u'' + x^2 u = \lambda u$ on the real axis.  Using Chebfun for
the quasimatrices, and approximating the real axis by $[-8,8]$, we
can compute eigenvalues based on an $\infty\times 40$ rectangular
Chebyshev spectral discretization with the code below,
just six lines long.  The first three eigenvalues come out as
$1.0000000008$, $3.0000000113$, and $5.0000005634$, and this
accuracy can be improved by increasing $n$.

{\small
\begin{verbatim}

        n = 40;
        L = chebop(@(x,u) -diff(u,2) + x^2*u,[-8,8]);
        G = chebpoly(0:n-1,[-8,8]);
        [Q,R] = qr(G);
        A = Q'*(L*G); C = R;
        lam = sort(eig(A,C))

\end{verbatim}
\par}

\noindent By adjusting a few of the commands we get a code for
the corresponding fully discrete computation with $100\times 40$
matrices, 100 being the default number of points in the {\tt
linspace} command.  Chebfun is still used in this code segment, but only because it
offers a convenient way to construct a matrix of sampled Chebyshev
polynomials scaled to $[-8,8]$ and their second derivatives.
The first three eigenvalues come out with approximately the
same accuracy as before as $1.0000000004$, $3.0000000050$,
and $5.0000002819$.

{\small
\begin{verbatim}

        n = 40;
        L = chebop(@(x,u) -diff(u,2) + x^2*u,[-8,8]);
        G = chebpoly(0:n-1,[-8,8]); LG = L*G;
        X = linspace(-8,8)'; G = G(X); LG = LG(X);
        [Q,R] = qr(G,0);
        A = Q'*LG; C = R;
        lam = sort(eig(A,C))

\end{verbatim}
\par}

For both of the computations just presented, the accuracy of the
computed eigenvalues is undiminished if the formulation
(\ref{unstable}) without the QR factorization is used instead of (\ref{new}).
This makes sense since $G$ is a matrix of Chebyshev polynomials on $[-1,1]$,
hence well-conditioned. 

\smallskip

{\em Example 2.  Wave oscillator with two boundary conditions.}
Our second example, mentioned on p.~\pageref{p3b}, is
$-(4/\pi^2)u'' = \lambda u$ on $[-1,1]$ with $u(\pm 1) = 0$.
The following code implements an $\infty\times 30$ quasimatrix
discretization, computing the first ten eigenvalues $1, 4, 9,\dots,
100$ to 11--14 digits of relative accuracy.

{\small
\begin{verbatim}

        n = 30;
        L = chebop(@(x,u) -(4/pi^2)*diff(u,2));
        G = chebpoly(0:n-1);
        [Q,R] = qr(G);
        A = [Q(:,1:n-2)'*(L*G); G(-1); G(1)];
        C = [R(1:n-2,:); zeros(2,n)];
        lam = sort(eig(A,C))

\end{verbatim}
\par}

\noindent Here is the adjustment needed for a fully discrete
$200\times 30$ discretization using equispaced points in
$[-1,1]$.  (It makes little difference if 200 Chebyshev points are
used instead, since we are in the regime $m\gg n$ of least-squares
sampling with plenty of sample points.  With $m\approx n$, it would be
important to be careful about the distribution of sample points, but
rectangular numerical methods make it unnecessary for $m$ to be small.)
The relative accuracy of the first ten eigenvalues
is now 9--14 digits, which returns to 11--14 digits if $n$ is
increased to $34$.

{\small
\begin{verbatim}

        n = 30;
        L = chebop(@(x,u) -(4/pi^2)*diff(u,2));
        G = chebpoly(0:n-1); LG = L*G;
        X = linspace(-1,1,200)'; G = G(X); LG = LG(X);
        [Q,R] = qr(G,0);
        A = [Q(:,1:n-2)'*LG; G(1,:); G(end,:)];
        C = [R(1:n-2,:); zeros(2,n)];
        lam = sort(eig(A,C))

\end{verbatim}
\par}

As with the last pair of computations, there is again little difference
in accuracy here if one bypasses the QR factorization and uses
(\ref{unstable}) instead of (\ref{new}).

These examples are of a Chebyshev spectral
flavor, and in such cases, at least in simple domains, 
square discretizations are often
readily available.  Now we turn to problems related to the Method
of Fundamental Solutions or RBF or other meshfree
discretizations, where the need for rectangular formulations
is more pressing.  The reason is that the representation of the
solution involves $n$ points that do not lie in the domain, hence
have no naturally associated grid for interpolation or quadrature.

\smallskip

{\em Example 3.  Wave oscillator, method of fundamental solutions.}
Our third example is the problem $-(4/\pi^2)u'' = \lambda u$
on $[-1,1]$ with $u(\pm 1) = 0$ again, but now solved by a kind of
method of fundamental solutions, with the solution represented as a linear
combination of point charges.  The following code implements
a $150 \times 35$ matrix discretization involving a constant
term plus 34 point charge potentials $\log|x-p_j|$ with $p_j$
equally spaced from $-1.5+0.5\kern .3pt i$ to $1.5+0.5\kern .3pt i$.  The first ten
eigenvalues $1, 4, 9,\dots, 100$ are computed to 12--15 digits
of relative accuracy.  A similar $150\times 35$ discretization
based on 17 dipoles equally spaced from $-1.5+0.5\kern .3pt i$ to $1.5+0.5\kern .3pt i$,
that is, real and imaginary parts of complex poles $1/(x-p_j)$
(not shown), gives 10--11 digits.

{\small
\begin{verbatim}

        n = 35;
        pts = linspace(-1.5+.5i,1.5+.5i,n-1);
        X = linspace(-1,1,150)';
        G = [X.^0 log(abs(X-pts))]; Gpp = [0*X -real(1./(X-pts).^2)];
        LG = -(4/pi^2)*Gpp;
        [Q,R] = qr(G,0); A = [Q(:,1:n-2)'*LG; G(1,:); G(end,:)];
        C = [R(1:n-2,:); zeros(2,n)];
        lam = sort(eig(A,C));

\end{verbatim}
\par}

For this problem, the QR factorization makes a big difference.  
If we use (\ref{unstable}) instead of (\ref{new}),
some spurious eigenvalues appear and the first ten nonspurious computed eigenvalues 
fall to 2--10 digits of accuracy.

\smallskip

{\em Example 4.  Quantum oscillator with singularity: lightning discretization.}
We now look at a problem with a singularity,
\begin{equation}
-0.01 u'' + |x|^{1/2} u = \lambda u, \quad u(\pm 1) = 0,
\label{sing}
\end{equation}
posed on the interval $[-1,1]$.  This is a Schr\"odinger equation with
the singular potential $V(x) = |x|^{1/2}$.  Smooth discretizations will have difficulty
achieving more than around 3 digits of accuracy, but we can do better with
a ``lightning discretization'' involving poles exponentially clustered
near the singular point $x=0$.  Specifically, following~eq.~(3.2) of~\cite{lightning},
a formula that is justified
in~\cite{clustering}, we fix a number $\npoles\ge 0$ of poles and define
\begin{equation}
d_j = \exp(4(\sqrt j -\sqrt{\vphantom{j}\npoles}\kern 1.5pt)), \quad 1 \le j\le \npoles.
\label{poles}
\end{equation}
The columns of $G$ will include both the real and the imaginary parts of
the simple pole functions $d_j^3/(x-i\kern .3pt d_j)$, making $2\kern .5pt\npoles$ columns all
together.  (The constant $d_j^3$ is included for scaling, since
the second derviative of this function is $2d_j^3/(x-id_j)^3$.)
In addition we fix a number $\npoly\ge 0$ and include
the Chebyshev polynomials $T_k(x)$ with $0 \le k \le \npoly$ as further columns of the 
matrix.

\begin{figure}
\begin{center}
\vskip 10pt
\includegraphics[scale=.75]{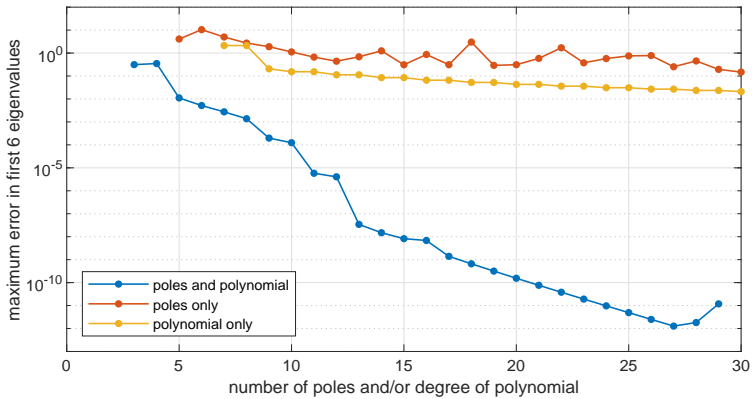}
\vskip 4pt
\end{center}
\caption{\label{singfig}Maximum error in the first six computed eigenvalues of the singular Schr\"odinger
problem~$(\ref{sing})$.  A rectangular discretization combining a smooth polyomial term with
exponentially clustered ``lightning'' poles converges rapidly to $11$ digits of accuracy
with a basis matrix $G$ with $4000$ rows and $3\times 25 + 1 = 76$ columns.
Neither the polynomial nor the clustered poles alone get better than $3$ digits.}
\end{figure}

Figure~\ref{singfig} shows results for this scheme for three sequences of computations in
which $\npoly$, or $\npoles$, or both range from $0$ to $30$.
The interval $[-1,1]$ is discretized by 3000 points exponentially spaced
from $10^{-10}$ to $1$ and their negatives, so
the matrices have $6000$ rows and between $6$ and $91$ columns.  (This space
discretization could undoubtedly be improved.)
Neither poles nor smooth polynomials alone give good accuracy, but in combination
they achieve up to 12 digits before a plateau is reached.

If Figure~\ref{singfig} is recomputed based on the formulation (\ref{unstable}) without QR factorization,
the accuracy falls to 
2 digits or so for $n>15$ and spurious eigenvalues appear.  For our further
examples we will report results only from the stable formulation (\ref{new}).

\section{\label{2D}Two-dimensional examples (PDE\kern .5pt s)}
Now we move to two-dimensional (2D) domains and PDE eigenvalue
problems.  Though problems without boundaries can certainly
be considered (such as the 2D harmonic oscillator in the $x$-$y$ plane), we shall look at examples
where $\Omega$ has a boundary $\dO$ with explicit boundary conditions, leading to rectangular
discretizations of the Variant~3 form (\ref{var3}).

\smallskip

{\em Example 5.  Circular drum, RBF discretization.}
Consider the planar Laplace problem mentioned on p.~\pageref{bessel},
\begin{equation}
-\Delta u = \lambda u, \quad |z|<1,
\label{drum}
\end{equation}
with $u(z) = 0$ for $|z|=1$.
The eigenvalues are the squares of the zeros of the Bessel functions $J_k(r)$, 
$k\ge 0$.  For $k=0$, the
eigenfunctions are axisymmetric and the eigenvalues are simple,
whereas for $k\ge 1$, the eigenfunctions are not axisymmetric and each
eigenvalue is of multiplicity~2.

\begin{figure}
\begin{center}
\includegraphics[scale=.52]{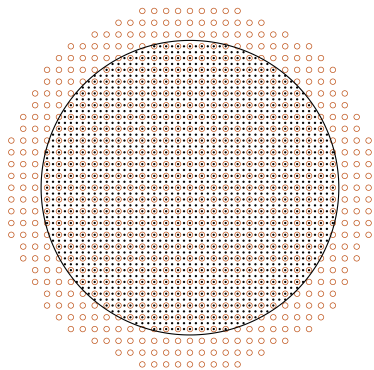}
\vskip 5pt
\end{center}
\caption{\label{RBFgrids}RBF discretization of Example~5.  The dots show $1941$ interior sample
points in the unit disk, and the circles show $769$ RBF centers in the disk $|z|\le 1.25$.
The unit circle boundary is discretized by $400$ equispaced points.  The
resulting rectangular eigenproblem is of dimensions $2241\times 770$.}
\end{figure}

\begin{figure}
\vskip -9pt
\indent\kern -18pt\includegraphics[scale=.85]{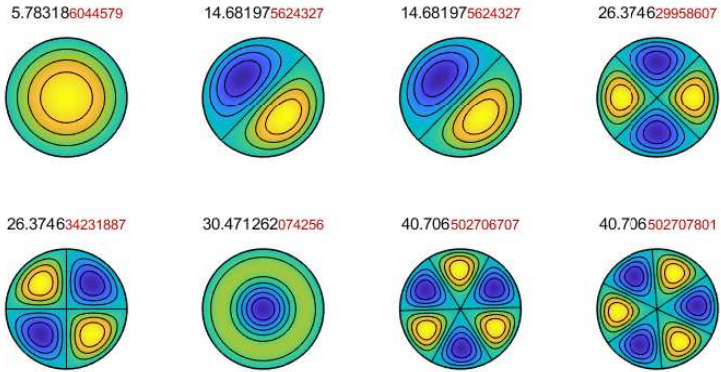}
\vspace{-114pt}
\caption{\label{RBFfig}First 8 computed eigenvalues and
eigenfunctions of a disk, based on an RBF discretization defined by
(\ref{multiq})--(\ref{RBFformula}) with the RBF and sampling grids of Figure~$\ref{RBFgrids}$.
The rectangular eigenvalue problem is of dimensions $2241\times 770$.
Correct digits are printed in black and incorrect ones in a smaller font in red.
The
eigenvalues of multiplicity~$2$ are identified correctly, but the associated eigenfunction
pairs do not come out orthogonal, reflecting the fact that the numerical method is
not self-adjoint.}
\end{figure}

Our first rectangular discretization will be based on RBF\kern .5pt s.
We follow~\cite{plattedrisc} and take as a radial basis function
the multiquadric
\begin{equation}
\phi(r) = \sqrt{c^2 + r^2}
\label{multiq}
\end{equation}
for a fixed parameter $c$, so that each eigenfunction is approximated by
a sum
\begin{equation}
u(z) = a_0 + \sum_{k=1}^{n-1} a_k \phi(|z-\zeta_k|),
\label{RBFformula}
\end{equation}
where $\{\zeta_k\}$ is a set of $n-1$ centers.
In~\cite{plattedrisc}, as in most RBF literature, 
the emphasis is on obtaining square discretizations based on interpolation.
This requires care in selecting the centers, which must
be clustered near the boundary to avoid a Runge phenomenon~\cite{platte11}.
In rectangular mode, however,  with $m\gg n$ sample points, one can be more relaxed.
To illustrate the method, Figure~\ref{RBFgrids} shows a square grid of sample points
in the unit disk $|z|\le 1$ with spacing $0.04$
as well as a sparser square grid of RBF centers in the disk
$|\zeta|\le 1.25$ with spacing $0.08$.
We take $c=0.4$ for the constant of (\ref{multiq}).
Figure~\ref{RBFfig} shows that the resulting $2241\times 770$ rectangular
eigenvalue problem gives 5--8 digit accuracy in the first eight eigenvalues.

\begin{figure}
\vskip 2pt
\indent\kern -18pt\includegraphics[scale=.85]{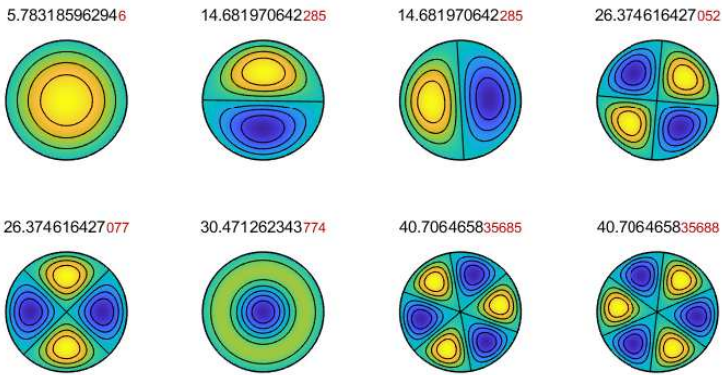}
\vspace{-114pt}
\caption{\label{diskfig}Like Figure~$\ref{RBFfig}$ but for a Fourier extension
discretization.
The rectangular eigenvalue problem is of dimensions $1545\times 221$.  The
degenerate pairs $2$--$3$, $4$--$5$, and $7$--$8$ again show orientations at
arbitrary angles.}
\end{figure}
\smallskip

{\em Example 6.  Circular drum, Fourier extension discretization.}
Consider (\ref{drum}) for a circular drum again,
but now discretized by a Fourier extension method.
For a rectangular discretization of (\ref{drum}), we start from a $41\times 41$ square grid in $[-1,1]^2$
(i.e., grid spacing $0.05$)
and discard the points outside the unit disk; the 1245 points that remain are
our interior sample points.  On the boundary we take 300 equispaced sample points.
The basis functions are the products
$\cos(kx)\cos(my)$, 
$\cos(kx)\sin(my)$, 
$\sin(kx)\cos(my)$, and
$\sin(kx)\sin(my)$ with $0\le k+m \le 10$ (discarding those that are exactly zero).
This leads to an eigenvalue problem of
dimensions $1545\times 221$, and
Figure~\ref{diskfig} shows the first eight computed eigenfunctions, with
eigenvalues accurate to $9$--$11$ digits. 

\smallskip

\begin{figure}
\indent\kern -18pt\includegraphics[scale=.88]{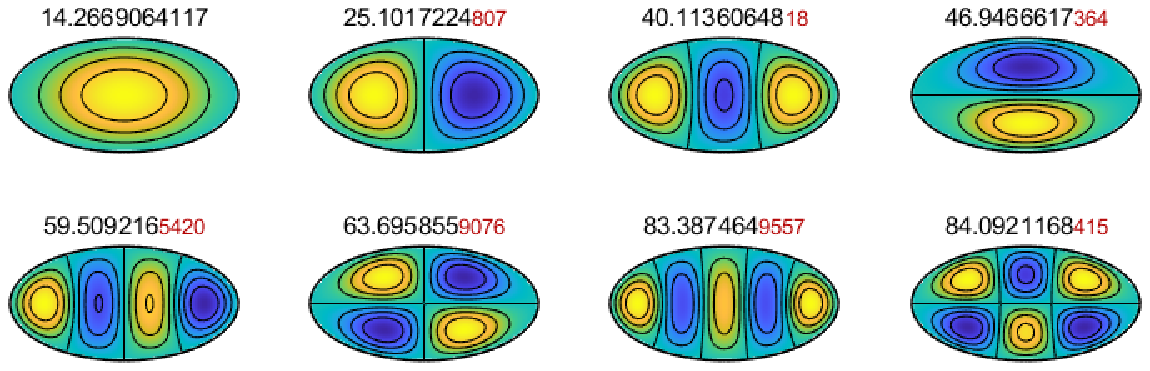}
\vspace{-186pt}
\caption{\label{ellipsefig}Like Figure~\ref{diskfig} but for an elliptical drum
of axis lengths $1$ and $\textstyle{1\over 2}$.}
\end{figure}

{\em Example 7.  Elliptical drum, Fourier extension discretization.}
Figure~\ref{ellipsefig} shows results for another 2D Fourier extension
computation, this time involving an elliptical drum of length $1$ and width $1/2$.
Although it doesn't make much difference for this problem, we have
switched here to a representation in which the 
basis functions are orthogonalized by a Vandermonde with Arnoldi
process~\cite[Example 3]{VA}.  Two-dimensional Vandermonde with Arnoldi has
been utilized previously for bivariate and trivariate
polynomials in~\cite{austin}, \cite{hokanson}, and~\cite{zhu},
and here we do it for Fourier extension.
(In separate experiments not reported here, we have successfully computed
eigenmodes of the ellipse in this manner by bivariate polynomials.)
Setting $X = e^{ix}$ and $Y = e^{iy}$, we note that 
the real part of $X^k Y^{\pm m}$ is $\cos(kx)\cos(my) \mp \sin(kx)\sin(my)$ and 
the imaginary part is $\pm \cos(kx)\sin(my) + \sin(kx)\cos(my)$, so
these real and imaginary parts span the necessary space of bivariate
trigonometric polynomials.
To be precise, we fix $K\ge 1$ and work with integers $k$ and $m$ with
$0\le k \le K$ and $0\le m \le K$ for $k=0$, 
$k-K \le m \le K-k$ for $k\ge 1$.  Arnoldi orthogonalization
is carried out in the order $1, Y, X, Y^2, XY, X^2, Y^3, XY^2, X^2Y, X^3,\dots$
(compare the paragraph after eq.~(7) of~\cite{austin}).
The rectangular matrix whose eigenfunctions are shown in the figure
is of dimensions $1399\times 313$.

\begin{figure}
\begin{center}
\includegraphics[scale=.72]{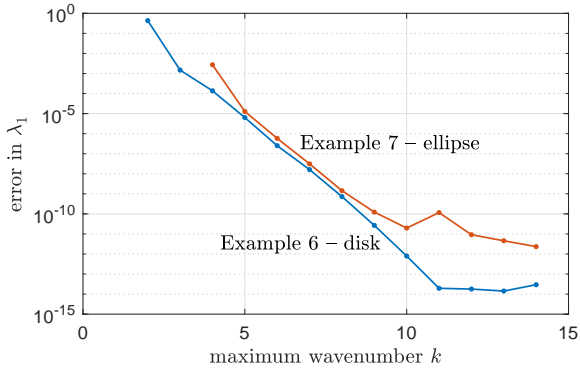}
\end{center}
\vskip 6pt
\caption{\label{figrates}Convergence curves for the two examples involving
Fourier extension discretizations.  Exponential convergence is observed to
around 14 digits for the disk and 11 digits for the ellipse.}
\end{figure}

Fourier approximations of analytic functions on analytic domains should
converge exponentially, and for Examples~6 and 7, the data confirm this
nicely, as shown in Figure~\ref{figrates}.
(Here and in the next example, the correct eigenvalues are not known analytically but
are determined numerically by higher-resolution calculations.)
For more complicated domains, however, especially if they are nonconvex,
the exponential rate becomes very slow because the solution can only be
analytically continued a short distance outside the boundary.  Such effects
have been studied by Barnett and Betcke~\cite{bb}, and for more on the
theory of analytic continuation of Helmholtz fields see~\cite{millar}.
This difficulty pertains to the choice of expansion functions,
not to the the method of dealing with them by rectangular eigenvalue problems.

\smallskip

{\em Example 8.  L-shaped region with singular terms.}
Solutions to PDE\kern .5pt s in regions with corners usually have
corner singularities, which make it challenging to get 
high accuracy.  In the context of the global representations
explored in this paper, a natural idea for such problems would be to 
combine a general purpose set of basis functions to capture the
``smooth part'' of the solution with additional singular terms near
the corners.   For Laplace Dirichlet or Neumann
problems, representations of this kind
led to the lightning and log-lightning solvers introduced in~\cite{lightning}
and~\cite{loglightning}.  Here we illustrate that such an approach may be
effective for eigenvalue problems too.  This is a PDE analogue of Example~4
for ODE\kern .5pt s.

\begin{figure}
\indent\kern -23pt\includegraphics[scale=.83]{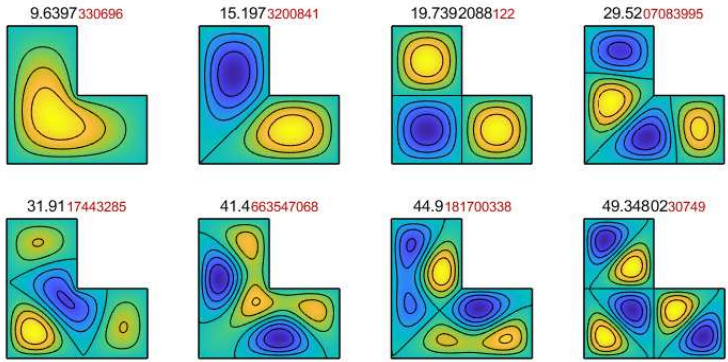}
\vspace{-112pt}
\caption{\label{Lfig}Eigenfunctions of an L-shaped region computed by a
$1617\times 343$ eigenvalue problem, with 30 of the matrix columns devoted
to resolving the singularity at the reentrant corner.}
\end{figure}

Our example, shown in Figure~\ref{Lfig}, is the planar drum
(\ref{drum}) in the form of the L-shaped region
well known from the MATLAB logo, the square $[-1,1]^2$ with one
corner removed.  (For numerical eigenvalues of
this and other drums calculated by a more
specialized method to an accuracy of 8 digits, see~\cite{computed}.)
The boundary has been discretized by 420 points exponentially
clustered near the reentrant corner, and the interior by a square grid of
spacing $1/20$.  This gives a $1617\times 313$ matrix, to which a further
30 columns are added corresponding to the corner singularity.
These 30 terms are chosen to capture the dominant behavior of the functions
$J_{(2/3)j}(\sqrt \lambda \kern .9pt r) \sin((2/3)ij\theta)$ that arise in
series expansions of eigenfunctions near reentrant right-angle corners, where
$r$ is the distance from the
corner and $\theta$ is the angle measured from one of the adjacent sides.
Specifically, we used discretizations of the 30 functions
$r^a \sin(b\kern .7pt \theta)$ with $a$ and $b$ given by
$$
 a = 2/3, 4/3, ~~ 8/3, 10/3, ~~ 14/3, 16/3, ~~ 20/3, 22/3, ~~ 26/3, 28/3
$$
and
$$
 b = a, \, a-2, \, a-4, \dots > 0.
$$
The figure shows that this rectangular discretization 
computes the first 8 eigenvalues to about 4 digits of accuracy.
So far as we know, discretizations of this kind have not been considered
before for eigenvalue problems, and we hope to present them more fully, and improve them,
in a future publication.

\section{\label{disc}Discussion}
The most robust discretizations of differential equation
eigenvalue problems, and the ones
with the strongest theoretical support, often involve
square matrices, especially in the self-adjoint case.  The theory of
finite element methods has brought such discretizations to an advanced state.

For some problems, however, whether because of irregular geometry, 
nonself-adjointness, or the presence of singularities, a good square
matrix discretization may not be readily available.
The aim of this paper has been to show that in such cases rectangular
matrices may offer an eminently practical alternative, often making
possible high accuracy solutions with a global representation (hence
perfectly smooth in the interior, and very fast to evaluate).
We make no claim of guaranteed success, and indeed, in most
of our experiments, which are based on new kinds of discretizations
with little previous literature, it has been necessary to try several parameter choices to 
get good accuracy and avoid spurious modes.  
With further work, more may be learned about these matters and rectangular eigenvalue
methods may be developed with guarantees of robustness and accuracy.
These methods are easy and flexible and deserve ongoing attention.

\bmhead{Acknowledgements}
We are grateful for helpful suggestions from Alex Barnett,
Timo Betcke, Toby
Driscoll, Mark Embree, Greg Fasshauer, Abi Gopal, Dave Hewett, Norm
Levenberg, Rodrigo Platte, Euan Spence, and Alex Townsend.

\section*{\label{coi}Declarations}

{\bf Conflicts of interest. }  The authors declare that they have no conflicts of interest.

\end{document}